\documentclass[12pt]{article}
\usepackage{amssymb,amsmath}

\newcommand{\Z}{\mathbb{Z}}

\newcommand{\F}{\mathbb{F}}
\newcommand{\OO}{\mathcal{O}}
\newcommand{\M}{\mathcal{M}}

\newcommand{\Aut}{\mathrm{Aut}}
\newcommand{\sgn}{\mathrm{sgn}}
\newcommand{\id}{\mathrm{id}}
\newcommand{\norm}{\mathrm{N}}
\newcommand{\Tr}{\mathrm{Tr}}
\newcommand{\Kbar}{\overline{K}}
\newcommand{\vbar}{\overline{v}}
\newcommand{\phibar}{\overline{\phi}}
\newcommand{\betabar}{\overline{\beta}}
\newcommand{\mubold}{\boldsymbol{\mu}}
\newcommand{\nubold}{\boldsymbol{\nu}}
\newcommand{\kappabold}{\boldsymbol{\kappa}}
\newcommand{\lambdabold}{\boldsymbol{\lambda}}
\newcommand{\proof}{\noindent{\em Proof: }}
\newcommand{\qed}{\hspace{\fill}$\square$}
\newcommand{\ra}{\rightarrow}

\newcommand{\dst}{\displaystyle}

\DeclareMathOperator{\lcm}{lcm}
\DeclareMathOperator{\Char}{char}

\newtheorem{theorem}{Theorem}
\newtheorem{lemma}[theorem]{Lemma}
\newtheorem{prop}[theorem]{Proposition}
\newtheorem{cor}[theorem]{Corollary}

\newenvironment{remark}{\refstepcounter{theorem}
\noindent{\bf Remark \thesection.\arabic{theorem}} }{}

\newenvironment{example}{\refstepcounter{theorem}
\noindent{\bf Example \thesection.\arabic{theorem}} }{}

\numberwithin{equation}{section}
\numberwithin{theorem}{section}

\title{Extensions of local fields and elementary
symmetric polynomials}

\author{Kevin Keating \\
Department of Mathematics \\
University of Florida \\
Gainesville, FL 32611 \\
USA \\[.2cm]
{\tt keating@ufl.edu}}

\begin{document}

\maketitle

\begin{abstract}
\noindent
Let $K$ be a local field whose residue field has
characteristic $p$ and let $L/K$ be a finite separable
totally ramified extension of degree $n=up^{\nu}$.  Let
$\sigma_1,\dots,\sigma_n$ denote the $K$-embeddings of
$L$ into a separable closure $K^{sep}$ of $K$.  For
$1\le h\le n$ let $e_h(X_1,\dots,X_n)$ denote the $h$th
elementary symmetric polynomial in $n$ variables, and
for $\alpha\in L$ set $E_h(\alpha)
=e_h(\sigma_1(\alpha),\dots,\sigma_n(\alpha))$.  Set
$j=\min\{v_p(h),\nu\}$.  We show that for $r\in\Z$ we
have $E_h(\M_L^r)\subset \M_K^{\lceil(i_j+hr)/n\rceil}$,
where $i_j$ is the $j$th index of inseparability of
$L/K$.  In certain cases we also show that $E_h(\M_L^r)$
is not contained in any higher power of $\M_K$.
\end{abstract}

\section{The problem} \label{prob}

Let $K$ be a field which is complete with respect to a
discrete valuation $v_K$.  Let $\OO_K$ be the ring of
integers of $K$ and let $\M_K$ be the maximal ideal of
$\OO_K$.  Assume that the residue field
$\Kbar=\OO_K/\M_K$ of $K$ is a perfect field of
characteristic $p$.  Let $K^{sep}$ be a separable
closure of $K$, and let $L/K$ be a finite totally
ramified subextension of $K^{sep}/K$ of degree
$n=up^{\nu}$, with $p\nmid u$.  Let
$\sigma_1,\dots,\sigma_n$ denote the $K$-embeddings of
$L$ into $K^{sep}$.  For $1\le h\le n$ let
$e_h(X_1,\dots,X_n)$ denote the $h$th elementary
symmetric polynomial in $n$ variables, and define
$E_h:L\ra K$ by setting $E_h(\alpha)
=e_h(\sigma_1(\alpha),\dots,\sigma_n(\alpha))$ for
$\alpha\in L$.  We are interested in the relation
between $v_L(\alpha)$ and $v_K(E_h(\alpha))$.  In
particular, for $r\in\Z$ we would like to compute the
value of
\[g_h(r)=\min\{v_K(E_h(\alpha)):\alpha\in\M_L^r\}.\]
The following proposition shows that $g_h(r)$ is a
well-defined integer:

\begin{prop}
Let $L/K$ be a totally ramified extension of degree $n$.
Let $r\in\Z$ and let $h$ satisfy $1\le h\le n$.  Then
$E_h(\M_L^r)\subset\M_K^{\lceil hr/n\rceil}$ and
$E_h(\M_L^r)\not=\{0\}$.
\end{prop}

\proof For the first claim we observe that if
$\alpha\in\M_L^r$ then $v_L(E_h(\alpha))\ge hr$, and
hence $v_K(E_h(\alpha))\ge hr/n$.
To prove the second claim let $A=L\otimes_KK^{sep}$ and
define $\tilde{E}_h:A\ra K^{sep}$ as follows.  For
$\beta\in A$ define $T_{\beta}:A\ra A$ by
$T_{\beta}(x)=\beta x$ for $x\in A$.  Then
$T_{\beta}$ is a $K^{sep}$-linear map.  Let
\[\det(X\cdot\id_A-T_{\beta})=X^n-c_1X^{n-1}+\dots
+(-1)^{n-1}c_{n-1}X+(-1)^nc_n\]
be the characteristic polynomial of $T_{\beta}$, and
set $\tilde{E}_h(\beta)=c_h$.  Since $L/K$ is separable
we have an isomorphism of $K^{sep}$-algebras
$A\cong(K^{sep})^n$.  It follows that $\tilde{E}_h$ is
onto.  There is an embedding of $K$-algebras $i:L\ra A$
defined by $i(\alpha)=\alpha\otimes1$ for $\alpha\in L$.
It follows from the definitions that
$\tilde{E}_h\circ i=E_h$.

     Let $\{v_1,\dots,v_n\}$ be a basis for $L$ over
$K$.  Then $\{v_1\otimes1,\dots,v_n\otimes1\}$ is a
basis for $A$ over $K^{sep}$.  For $x_i\in K^{sep}$
define
\[F(x_1,\dots,x_n)
=\tilde{E}_h(v_1\otimes x_1+\dots+v_n\otimes x_n).\]
Then $F$ is a degree-$h$ form on $(K^{sep})^n$.
Furthermore, since $\tilde{E}_h$ is onto, $F$ in
nontrivial.  Since $K$ is an infinite field there are
$d_i\in K$ such that $F(d_1,\dots,d_n)\not=0$.  Set
$\alpha=d_1v_1+\dots+d_nv_n$.  Then $\alpha\in L$ and
\[E_h(\alpha)=\tilde{E}_h(\alpha\otimes1)=F(d_1,\dots,d_n)
\not=0.\]
Let $\pi_K$ be a uniformizer for $K$.  Then for $t$
sufficiently large we have $\pi_K^t\alpha\in\M_L^r$ and
$E_h(\pi_K^t\alpha)=\pi_K^{ht}E_h(\alpha)\not=0$.  Hence
$E_h(\M_L^r)\not=\{0\}$. \qed \medskip

     Since $L/K$ is totally ramified we have
\[v_K(E_n(\alpha))=v_K(\norm_{L/K}(\alpha))=v_L(\alpha),\]
and hence $g_n(r)=r$ for $r\in\Z$.  The map
$E_1=\Tr_{L/K}$ is also well-understood, at least when
$L/K$ is a Galois extension of degree $p$ (see
\cite[V\,\S3, Lemma 4]{cl} or \cite[III,
Prop.\,1.4]{FV}).

\begin{prop} \label{different}
Let $L/K$ be a totally ramified extension of degree $n$
and let $\M_L^d$ be the different of $L/K$.  Then for
every $r\in\Z$ we have
$E_1(\M_L^r)=\M_K^{\lfloor(d+r)/n\rfloor}$.  Therefore
$g_1(r)=\lfloor(d+r)/n\rfloor$.
\end{prop}

\proof Since $E_1(\M_L^r)$ is a nonzero fractional ideal
of $K$ we have $E_1(\M_L^r)=\M_K^s$ for some $z\in\Z$.
By Proposition~7 in \cite[III\,\S3]{cl} we have
\begin{align*}
\M_L^{d+r}&\subset\OO_L\cdot\M_K^s=\M_L^{ns} \\
\M_L^{d+r}&\not\subset\OO_L\cdot\M_K^{s+1}=\M_L^{n(s+1)}.
\end{align*}
It follows that $d+r\ge ns$ and $d+r<n(s+1)$, and hence
that $s=\lfloor(d+r)/n\rfloor$.~\qed \medskip

     In this paper we determine a lower bound for
$g_h(r)$ which depends on the indices of inseparability
of $L/K$.  When $h=p^j$ with $0\le j\le\nu$ and $\Kbar$
is large enough we show that $g_h(r)$ is equal to this
lower bound.  This leads to a formula for $g_{p^j}(r)$
which can be expressed in terms of a generalization of
the different of $L/K$ (see Remark~\ref{dj}).

     In sections \ref{digraphs} and \ref{subrings} we
prove some preliminary results involving symmetric
polynomials.  The main focus is on expressing monomial
symmetric polynomials in terms of elementary symmetric
polynomials.  In section~\ref{containment} we prove our
lower bound for $g_h(r)$.  In section~\ref{equal} we
show that $g_h(r)$ is equal to this lower bound in some
special cases.

\section{Symmetric polynomials and cycle digraphs}
\label{digraphs}

Let $n\ge1$, let $w\ge1$, and let $\lambdabold$ be a
partition of $w$.  We view $\lambdabold$ as a multiset
of positive integers
such that the sum $\Sigma(\lambdabold)$ of the elements
of $\lambdabold$ is equal to $w$.  The cardinality of
$\lambdabold$ is denoted by $|\lambdabold|$.  For
$k\ge1$ we let $k*\lambdabold$ be the partition of $kw$
which is the multiset sum of $k$ copies of
$\lambdabold$, and we let $k\cdot\lambdabold$ be the
partition of $kw$ obtained by multiplying the parts of
$\lambdabold$ by $k$.  If $|\lambdabold|\le n$ let
$m_{\lambdabold}(X_1,\dots,X_n)$ be the monomial
symmetric polynomial in $n$ variables associated to
$\lambdabold$.  For $1\le h\le n$ let
$e_h(X_1,\dots,X_n)$ denote the $h$th elementary
symmetric polynomial in $n$ variables.

     Let $r\ge1$ and let
$\phi(X)=a_rX^r+a_{r+1}X^{r+1}+\cdots$ be a power series
with generic coefficients $a_i$.  For a partition
$\mubold=\{\mu_1,\dots,\mu_h\}$ whose parts satisfy
$\mu_i\ge r$ set $a_{\mubold}=a_{\mu_1}a_{\mu_2}\dots
a_{\mu_h}$.  Then for $1\le h\le n$ we have
\begin{equation} \label{eh}
e_h(\phi(X_1),\dots,\phi(X_n))=\sum_{\mubold}
a_{\mubold}m_{\mubold}(X_1,\dots,X_n),
\end{equation}
where the sum ranges over all partitions $\mubold$ with
$h$ parts, all of which are $\ge r$.  By the fundamental
theorem of symmetric polynomials there is
$\psi_{\mubold}\in\Z[X_1,\dots,X_n]$ such that
$m_{\mubold}=\psi_{\mubold}(e_1,\dots,e_n)$.  In this
section we use a theorem of Kulikauskas and Remmel
\cite{kr} to compute certain coefficients of the
polynomials $\psi_{\mubold}$.

     The formula of Kulikauskas and Remmel can be
expressed in terms of tilings of a certain type of
digraph.  We say that a directed graph $\Gamma$ is a
cycle digraph if it is a disjoint union of finitely many
directed cycles of length $\ge1$.  We denote the vertex
set of $\Gamma$ by $V(\Gamma)$, and we define the sign
of $\Gamma$ to be $\sgn(\Gamma)=(-1)^{w-c}$, where
$w=|V(\Gamma)|$ and $c$ is the number of cycles that
make up $\Gamma$.

     Let $\Gamma$ be a cycle digraph with $w\ge1$
vertices and let $\lambdabold$ be a partition of $w$.
A $\lambdabold$-tiling of $\Gamma$ is a set $S$ of
subgraphs of $\Gamma$ such that
\begin{enumerate}
\item Each $\gamma\in S$ is a directed path of length
$\ge0$.
\item The collection $\{V(\gamma):\gamma\in S\}$ forms a
partition of the set $V(\Gamma)$.
\item The multiset $\{|V(\gamma)|:\gamma\in S\}$ is
equal to $\lambdabold$.
\end{enumerate}
Let $\mubold$ be another partition of $w$.  A
$(\lambdabold,\mubold)$-tiling of $\Gamma$ is an ordered
pair $(S,T)$, where $S$ is a $\lambdabold$-tiling of
$\Gamma$ and $T$ is a $\mubold$-tiling of $\Gamma$.  Let
$\Gamma'$ be another cycle digraph with $w$ vertices and
let $(S',T')$ be a $(\lambdabold,\mubold)$-tiling of
$\Gamma'$.  An isomorphism from $(\Gamma,S,T)$ to
$(\Gamma',S',T')$ is an isomorphism of digraphs
$\theta:\Gamma\ra\Gamma'$ which carries $S$ onto $S'$
and $T$ onto $T'$.  Say that $(S,T)$ is an admissible
$(\lambdabold,\mubold)$-tiling of $\Gamma$ if
$(\Gamma,S,T)$ has no nontrivial automorphisms.  Say
that the $(\lambdabold,\mubold)$-tilings $(S,T)$ and
$(S',T')$ of $\Gamma$ are isomorphic if there exists an
isomorphism from $(\Gamma,S,T)$ to $(\Gamma,S',T')$.
Let $\eta_{\lambdabold\mubold}(\Gamma)$
denote the number of isomorphism classes of admissible
$(\lambdabold,\mubold)$-tilings of $\Gamma$.

     Let $w\ge1$ and let $\lambdabold,\mubold$ be
partitions of $w$.  Set
\begin{equation} \label{dlm}
d_{\lambdabold\mubold}
=(-1)^{|\lambdabold|+|\mubold|}\cdot\sum_{\Gamma}\,
\sgn(\Gamma)\eta_{\lambdabold\mubold}(\Gamma),
\end{equation}
where the sum is over all isomorphism classes of cycle
digraphs $\Gamma$ with $w$ vertices.  Since 
$\eta_{\mubold\lambdabold}=\eta_{\lambdabold\mubold}$ we
have $d_{\mubold\lambdabold}=d_{\lambdabold\mubold}$.
Kulikauskas and Remmel \cite[Th.\,1(ii)]{kr} proved the
following:

\begin{theorem} \label{krtheorem}
Let $n\ge1$, let $w\ge1$, and let $\mubold$ be a
partition of $w$ with at most $n$ parts.  Let
$\psi_{\mubold}$ be the unique element of
$\Z[X_1,\dots,X_n]$ such that
$m_{\mubold}=\psi_{\mubold}(e_1,\dots,e_n)$.  Then
\[\psi_{\mubold}(X_1,\dots,X_n)=\sum_{\lambdabold}
d_{\lambdabold\mubold}\cdot
X_{\lambda_1}X_{\lambda_2}\dots X_{\lambda_k},\]
where the sum is over all partitions
$\lambdabold=\{\lambda_1,\dots,\lambda_k\}$ of $w$ such
that $\lambda_i\le n$ for $1\le i\le k$.
\end{theorem}

     The remainder of this section is devoted to
computing the values of
$\eta_{\lambdabold\mubold}(\Gamma)$ and
$d_{\lambdabold\mubold}$ in some special cases.

\begin{prop} \label{cycle4}
Let $w\ge1$, let $\lambdabold,\mubold$ be partitions
of $w$, and let $\Gamma$ be a directed cycle of length
$w$.  Assume that $\Gamma$ has a $\lambdabold$-tiling
$S$ which is unique up to isomorphism, and that
$\Aut(\Gamma,S)$ is trivial.  Similarly, assume that
$\Gamma$ has a  $\mubold$-tiling $T$ which is unique up
to isomorphism, and that $\Aut(\Gamma,T)$ is trivial.
Then $\eta_{\lambdabold\mubold}(\Gamma)=w$.
\end{prop}

\proof For $0\le i<w$ let $S_i$ be the rotation of $S$
by $i$ steps.  Then the isomorphism classes of
$(\lambdabold,\mubold)$-tilings of $\Gamma$ are
represented by $(S_i,T)$ for $0\le i<w$.  Since
$\Aut(\Gamma,T)$ is trivial, all these tilings are
admissible. \qed

\begin{prop} \label{cycle0}
Let $a,b,c,\ell,m,w$ be positive integers such that
$\ell a=mb+c=w$ and $b\not=c$.  Let $\lambdabold$ be the
partition of $w$ consisting of $\ell$ copies of $a$,
let $\mubold$ be the partition of $w$ consisting of
$m$ copies of $b$ and 1 copy of $c$, and let $\Gamma$ be
a directed cycle of length $w$.  Then
$\eta_{\lambdabold\mubold}(\Gamma)=a$.
\end{prop}

\proof The cycle digraph $\Gamma$ has a
$\lambdabold$-tiling $S$ which us unique up to
isomorphism, and a $\mubold$-tiling $T$ which is unique
up to isomorphism.  For $0\le i<a$ let $S_i$ be the
rotation of $S$ by $i$ steps.  Then the isomorphism
classes of $(\lambdabold,\mubold)$-tilings of $\Gamma$
are represented by $(S_i,T)$ for $0\le i<a$.  Since
$\Aut(\Gamma,T)$ is trivial, all these tilings are
admissible. \qed

\begin{prop} \label{cycle1}
Let $b,c,m,w$ be positive integers such that $mb+c=w$
and $b\not=c$.  Let $\lambdabold$ be the partition of
$w$ consisting of 1 copy of $w$ and let $\mubold$ be the
partition of $w$ consisting of $m$ copies of $b$ and 1
copy of $c$.  Then $\dst d_{\lambdabold\mubold}
=(-1)^{w+m+1}w$.
\end{prop}

\proof If the cycle digraph $\Gamma$ has a
$\lambdabold$-tiling then $\Gamma$ consists of a single
cycle of length $w$.  Hence by (\ref{dlm}) we get
$d_{\lambdabold\mubold}=(-1)^{w+m+1}
\eta_{\lambdabold\mubold}(\Gamma)$.  It follows from
Proposition~\ref{cycle0} that
$\eta_{\lambdabold\mubold}(\Gamma)=w$.  Therefore
$d_{\lambdabold\mubold}=(-1)^{w+m+1}w$. \qed

\begin{prop} \label{cycle2}
Let $a,b,\ell,m,w$ be positive integers such that
$\ell a=mb=w$.  Let $\lambdabold$ be the partition of
$w$ consisting of $\ell$ copies of $a$, let $\mubold$ be
the partition of $w$ consisting of $m$ copies of $b$,
and let $\Gamma$ be a directed cycle of length $w$.
\\[.2cm]
(a) The number of isomorphism classes of
$(\lambdabold,\mubold)$-tilings of $\Gamma$ is
$\gcd(a,b)$. \\[.2cm]
(b) Let $(S,T)$ be a $(\lambdabold,\mubold)$-tiling of
$\Gamma$.  Then the order of $\Aut(\Gamma,S,T)$ is
$\gcd(\ell,m)$.
\end{prop}

\proof (a) Identify $V(\Gamma)$ with $\Z/w\Z$ and
consider the translation action of $b\Z/w\Z$ on
$(\Z/w\Z)/(a\Z/w\Z)$.  The isomorphism
classes of $(\lambdabold,\mubold)$-tilings of $\Gamma$
correspond to the orbits of this action, and these
orbits correspond to cosets of
$a\Z+b\Z=\gcd(a,b)\cdot\Z$ in $\Z$. \\[.2cm]
(b) The automorphisms of $(\Gamma,S,T)$ are rotations of
$\Gamma$ by $m$ steps, where $m$ is a multiple of both
$a$ and $b$.  Hence the number of automorphisms is
$w/\lcm(a,b)$, which is easily seen to be equal to
$\gcd(\ell,m)$. \qed \medskip

     The following proposition generalizes the second
part of \cite[Th.\,6]{kr}.

\begin{prop} \label{cycle3}
Let $a,b,\ell,m,w$ be positive integers such that
$\ell a=mb=w$.  Let $\lambdabold$ be the partition of $w$
consisting of $\ell$ copies of $a$ and let $\mubold$ be
the partition of $w$ consisting of $m$ copies of $b$.
Set $u=\gcd(a,b)$ and $v=\gcd(\ell,m)$.  Then $\dst
d_{\lambdabold\mubold}=(-1)^{w-v+\ell+m}\binom{u}{v}$.
In particular, if $u<v$ then $d_{\lambdabold\mubold}=0$.
\end{prop}

\proof Set $i=a/u$ and $j=b/u$.  Then $m=vi$ and
$\ell=vj$.  Let $\Gamma$ be a cycle digraph which
has an admissible $(\lambdabold,\mubold)$-tiling,
and let $\Gamma_0$ be one of the cycles which makes up
$\Gamma$.  Then the length of $\Gamma_0$ is divisible by
$\lcm(a,b)=uij$.  Suppose $\Gamma_0$ has length $k\cdot
uij$.  Let $\lambdabold_0$ be the partition of $kuij$
consisting of $kj$ copies of $a=ui$, and let $\mubold_0$
be the partition of $kuij$ consisting of $ki$ copies of
$b=uj$.  Then by Proposition~\ref{cycle2}(b) every
$(\lambdabold_0,\mubold_0)$-tiling of $\Gamma_0$ has
automorphism group of order $\gcd(ki,kj)=k$.  Since
$\Gamma$ has an admissible
$(\lambdabold,\mubold)$-tiling we must have $k=1$.
Therefore $\Gamma$ consists of $v$ cycles, each of
length $uij$.  By Proposition~\ref{cycle2}(a) the number
of isomorphism classes of
$(\lambdabold_0,\mubold_0)$-tilings of a $uij$-cycle
$\Gamma_0$ is $\gcd(a,b)=u$.  An admissible
$(\lambdabold,\mubold)$-tiling of $\Gamma$ consists of
$v$ nonisomorphic $(\lambdabold_0,\mubold_0)$-tilings of
$uij$-cycles.  Hence the number of isomorphism classes
of admissible $(\lambdabold,\mubold)$-tilings of
$\Gamma$ is $\dst\eta_{\lambdabold\mubold}(\Gamma)
=\binom{u}{v}$.  Hence by (\ref{dlm}) we get
$\dst d_{\lambdabold\mubold}
=(-1)^{w-v+\ell+m}\binom{u}{v}$. \qed

\section{Some subrings of $\Z[X_1,\dots,X_n]$}
\label{subrings}

Let $n\ge1$.  In some cases we can get information about
the coefficients $d_{\lambdabold\mubold}$ which appear
in the formula for $\psi_{\mubold}$ given in
Theorem~\ref{krtheorem} by working directly with the
ring $\Z[X_1,\dots,X_n]$.  In this section we define a
family of subrings of $\Z[X_1,\dots,X_n]$.  We then
study the $p$-adic properties of the coefficients
$d_{\lambdabold\mubold}$ by showing that for certain
partitions $\mubold$ the polynomial $\psi_{\mubold}$ is
an element of one of these subrings.

     For $k\ge0$ define a subring $R_k$ of
$\Z[X_1,\dots,X_n]$ by
\[R_k=\Z[X_1^{p^k},\dots,X_n^{p^k}]
+p\Z[X_1^{p^{k-1}},\dots,X_n^{p^{k-1}}]+\dots
+p^k\Z[X_1,\dots,X_n].\]
We can characterize $R_k$ as the set of
$F\in\Z[X_1,\dots,X_n]$ such that for $1\le i\le k$
there exists $F_i\in\Z[X_1,\dots,X_n]$ such that
\begin{equation} \label{fcong}
F(X_1,\dots,X_n)\equiv F_i(X_1^{p^i},\dots,X_n^{p^i})
\pmod{p^{k+1-i}}.
\end{equation}

\begin{lemma} \label{pl}
Let $k,\ell\ge0$ and let $F\in R_k$.  Then
$p^\ell F\in R_{k+\ell}$ and $F^{p^\ell}\in R_{k+\ell}$.
\end{lemma}

\proof The first claim is clear.  To prove the second
claim with $\ell=1$ we note that for $1\le i\le k$ it
follows from (\ref{fcong}) that
\[F(X_1,\dots,X_n)^p\equiv
F_i(X_1^{p^i},\dots,X_n^{p^i})^p\pmod{p^{k+2-i}}.\]
In particular, the case $i=k$ gives
\begin{alignat*}{2}
F(X_1,\dots,X_n)^p&\equiv
F_k(X_1^{p^k},\dots,X_n^{p^k})^p&&\pmod{p^2} \\
&\equiv F_k(X_1^{p^{k+1}},\dots,X_n^{p^{k+1}})&&\pmod{p}.
\end{alignat*}
It follows that $F^p\in R_{k+1}$.  By induction we get
$F^{p^\ell}\in R_{k+\ell}$ for $\ell\ge0$. \qed

\begin{lemma} \label{Rn1}
Let $k,\ell\ge0$ and let $F\in R_k$.  Then for any
$\psi_1,\dots,\psi_n\in R_\ell$ we have
$F(\psi_1,\dots,\psi_n)\in R_{k+\ell}$.
\end{lemma}

\proof Since $F\in R_k$ we have
\[F(X_1,\dots,X_n)=
\sum_{i=0}^kp^{k-i}\phi_i(X_1^{p^i},\dots,X_n^{p^i})\]
for some $\phi_i\in\Z[X_1,\dots,X_n]$.  Since
$\psi_j\in R_\ell$, by Lemma~\ref{pl} we get
$\psi_j^{p^i}\in R_{i+\ell}$.  Since $R_{i+\ell}$ is a
subring of $\Z[X_1,\dots,X_n]$ it follows that
$\phi_i(\psi_1^{p^i},\dots,\psi_n^{p^i})\in
R_{i+\ell}$.  By Lemma~\ref{pl} we get
$p^{k-i}\phi_i(\psi_1^{p^i},\dots,\psi_n^{p^i})\in
R_{k+\ell}$.  We conclude that
$F(\psi_1,\dots,\psi_n)\in R_{k+\ell}$. \qed

\begin{prop} \label{power}
Let $w\ge1$ and let $\lambdabold$ be a partition of $w$
with at most $n$ parts.  For $j\ge0$ let
$\lambdabold^j=p^j\cdot\lambdabold$.  Then
$\psi_{\lambdabold^j}\in R_j$.
\end{prop}

\proof We use induction on $j$.  The case $j=0$ is
trivial.  Let $j\ge0$ and assume that
$\psi_{\lambdabold^j}\in R_j$.  Since
$\lambdabold^{j+1}=p\cdot\lambdabold^j$ we get
\begin{align*}
m_{\lambdabold^{j+1}}(X_1,\dots,X_n)
&=m_{\lambdabold^j}(X_1^p,\dots,X_n^p) \\
&=\psi_{\lambdabold^j}(e_1(X_1^p,\dots,X_n^p),\dots,
e_n(X_1^p,\dots,X_n^p)).
\end{align*}
Since $X_j^p\in R_1$ it follows from Lemma~\ref{Rn1}
that
\[e_i(X_1^p,\dots,X_n^p)=\theta_i(e_1,\dots,e_n)\]
for some $\theta_i\in R_1$.  Therefore
\[\psi_{\lambdabold^{j+1}}(X_1,\dots,X_n)
=\psi_{\lambdabold^j}(\theta_1(X_1,\dots,X_n),
\dots,\theta_n(X_1,\dots,X_n)).\]
By Lemma~\ref{Rn1} we get
$\psi_{\lambdabold^{j+1}}\in R_{j+1}$. \qed

\begin{cor} \label{val}
Let $t\ge j\ge0$, let $w'\ge1$, and set $w=w'p^t$.
Let $\lambdabold'$ be a partition of $w'$ and set
$\lambdabold=p^t\cdot\lambdabold'$.  Let $\mubold$
be a partition of $w$ such that there does not exist a
partition $\mubold'$ with $\mubold=p^{j+1}*\mubold'$.
Then $p^{t-j}$ divides $d_{\lambdabold\mubold}$.  This
holds in particular if $p^{j+1}\nmid|\mubold|$.
\end{cor}

\proof Since $d_{\lambdabold\mubold}$ does not depend
on $n$ we may assume without loss of generality that
$n\ge w$.  It follows from this assumption that
$|\lambdabold|\le n$, so by Proposition~\ref{power} we
have $\psi_{\lambdabold}\in R_t$.  Since $w\le n$
the parts of $\mubold=\{\mu_1,\dots,\mu_h\}$ satisfy
$\mu_i\le n$ for $1\le i\le h$.  Therefore the
formula for $\psi_{\lambdabold}$ given by
Theorem~\ref{krtheorem} includes
the term $d_{\mubold\lambdabold}X_{\mu_1}X_{\mu_2}\dots
X_{\mu_h}$.  The assumption on $\mubold$ implies that
$X_{\mu_1}X_{\mu_2}\dots X_{\mu_h}$ is not a $p^{j+1}$
power.  Since $\psi_{\lambdabold}\in R_t$ this implies
that $p^{t-j}$ divides $d_{\mubold\lambdabold}$.  Since
$d_{\lambdabold\mubold}=d_{\mubold\lambdabold}$ we get
$p^{t-j}\mid d_{\lambdabold\mubold}$.~\qed

\begin{prop} \label{cong}
Let $w'\ge1$, $j\ge1$, and $t\ge0$.  Let
$\lambdabold'$, $\mubold'$ be partitions of $w'$ such
that the parts of $\lambdabold'$ are all divisible by
$p^t$.  Set $w=w'p^j$, so that
$\lambdabold=p^j\cdot\lambdabold'$ and
$\mubold=p^j*\mubold'$ are partitions of $w$.  Then
$d_{\lambdabold\mubold} \equiv d_{\lambdabold'\mubold'}
\pmod{p^{t+1}}$.
\end{prop}

\proof As in the proof of Corollary~\ref{val} we may
assume without loss of generality that $n\ge w'$.  Then
$|\lambdabold'|=|\lambdabold|\le n$.  It follows from
Proposition~\ref{power} that
$m_{\lambdabold'}=\psi_{\lambdabold'}(e_1,\dots,e_n)$
for some $\psi_{\lambdabold'}\in R_t$.  Using induction
on $k$ we see that for $1\le i\le n$ and $k\ge0$ we have
\[e_i(X_1^{p^j},\dots,X_n^{p^j})^{p^k}\equiv
e_i(X_1,\dots,X_n)^{p^{j+k}}\pmod{p^{k+1}}.\]
Since $\psi_{\lambdabold'}\in R_t$ it follows that
\begin{alignat*}{2}
m_{\lambdabold}(X_1,\dots,X_n)&=
m_{\lambdabold'}(X_1^{p^j},\dots,X_n^{p^j}) \\
&=\psi_{\lambdabold'}(e_1(X_1^{p^j},\dots,X_n^{p^j}),\dots,
e_n(X_1^{p^j},\dots,X_n^{p^j})) \\
&\equiv\psi_{\lambdabold'}(e_1(X_1,\dots,X_n)^{p^j},\dots,
e_n(X_1,\dots,X_n)^{p^j})&&\pmod{p^{t+1}}.
\end{alignat*}
We also have $m_{\lambdabold}
=\psi_{\lambdabold}(e_1,\dots,e_n)$.  Hence by the
fundamental theorem of symmetric polynomials we get
\[\psi_{\lambdabold}(X_1,\dots,X_n)\equiv
\psi_{\lambdabold'}(X_1^{p^j},\dots,X_n^{p^j})
\pmod{p^{t+1}}.\]
Since $w'\le n$ the parts of $\mubold'$ and $\mubold$
are all $\le n$.  Therefore the formula for
$\psi_{\lambdabold'}$ given by Theorem~\ref{krtheorem}
includes the term $d_{\mubold'\lambdabold'}
X_{\mu_1'}X_{\mu_2'}\dots X_{\mu_h'}$, and the formula
for $\psi_{\lambdabold}$ includes the term
\[d_{\mubold\lambdabold}X_{\mu_1}X_{\mu_2}\dots
X_{\mu_{p^jh}}=d_{\mubold\lambdabold}
X_{\mu_1'}^{p^j}X_{\mu_2'}^{p^j}\dots X_{\mu_h'}^{p^j}.\]
It follows that $d_{\mubold\lambdabold}\equiv
d_{\mubold'\lambdabold'}\pmod{p^{t+1}}$, and hence that
$d_{\lambdabold\mubold}\equiv d_{\lambdabold'\mubold'}
\pmod{p^{t+1}}$.~\qed

\section{Containment} \label{containment}

Let $L/K$ be a totally ramified extension of degree
$n=up^{\nu}$, with $p\nmid u$.  Let
$\sigma_1,\dots,\sigma_n$ be the $K$-embeddings of $L$
into $K^{sep}$.  Let $1\le h\le n$ and recall that
$E_h:L\ra K$ is defined by
$E_h(\alpha)=e_h(\sigma_1(\alpha),\dots,\sigma_n(\alpha))$
for $\alpha\in L$.  In this section we define a function
$\gamma_h:\Z\ra\Z$ such that for $r\in\Z$ we have
$E_h(\M_L^r)\subset\M_K^{\gamma_h(r)}$.  The function
$\gamma_h$ will be defined in terms of the indices of
inseparability of the extension $L/K$.  In the next
section we show that $\OO_K\cdot E_h(\M_L^r)
=\M_K^{\gamma_h(r)}$ holds in certain cases.

     Let $\pi_L$ be a uniformizer for $L$ and let
\[f(X)=X^n-c_1X^{n-1}+\cdots+(-1)^{n-1}c_{n-1}X+(-1)^nc_n\]
be the minimum polynomial of $\pi_L$ over $K$.  Then
$c_h=E_h(\pi_L)$.  For $k\in\Z$ define
$\vbar_p(k)=\min\{v_p(k),\nu\}$.  For $0\le j\le\nu$ set
\begin{align*}
i_j^{\pi_L}
&=\min\{nv_K(c_h)-h:1\le h\le n,\;\vbar_p(h)\le j\} \\
&=\min\{v_L(c_h\pi_L^{n-h}):
1\le h\le n,\;\vbar_p(h)\le j\}-n.
\end{align*}
Then $i_j^{\pi_L}$ is either a nonnegative integer or
$\infty$.  If $\Char(K)=p$ then $i_j^{\pi_L}$ must be
finite, since $L/K$ is separable.  If $i_j^{\pi_L}$ is
finite write $i_j^{\pi_L}=a_jn-b_j$ with
$1\le b_j\le n$.  Then $v_K(c_{b_j})=a_j$,
$v_K(c_h)\ge a_j$ for all $h$ with $1\le h<b_j$ and
$\vbar_p(h)\le j$, and $v_K(c_h)\ge a_j+1$ for all $h$
with $b_j<h\le n$ and $\vbar_p(h)\le j$.  Let
$e_L=v_L(p)$ denote the absolute ramification index of
$L$.  We define the $j$th index of inseparability of
$L/K$ to be
\begin{equation*}
i_j=\min\{i_{j'}^{\pi_L}+(j'-j)e_L:j\le j'\le\nu\}.
\end{equation*}
By Proposition~3.12 and Theorem~7.1 of \cite{heier},
$i_j$ does not depend on the choice of $\pi_L$.
Furthermore, our definition of $i_j$ agrees with
Definition~7.3 in \cite{heier} (see also
\cite[Remark~2.5]{towers}; for the characteristic-$p$
case see \cite[pp.\,232--233]{fried} and
\cite[\S2]{fm}).

     The following facts are easy consequences of the
definitions:
\begin{enumerate}
\item $0=i_{\nu}<i_{\nu-1}\le\dots\le i_1\le i_0<\infty$.
\item If $\Char(K)=p$ then $e_L=\infty$, and hence
$i_j=i_j^{\pi_L}$.
\item Let $m=\vbar_p(i_j)$.  If $m\le j$ then
$i_j=i_m=i_j^{\pi_L}=i_m^{\pi_L}$.  If $m>j$ then
$\Char(K)=0$ and $i_j=i_m^{\pi_L}+(m-j)e_L$.
\end{enumerate}

\begin{lemma} \label{bound}
Let $1\le h\le n$ and set $j=\vbar_p(h)$.  Then
$v_L(c_h)\ge i_j^{\pi_L}+h$, with equality if and only
if either $i_j^{\pi_L}=\infty$ or $i_j^{\pi_L}<\infty$
and $h=b_j$.
\end{lemma}

\proof If $i_j^{\pi_L}=\infty$ then we certainly have
$v_L(c_h)=\infty$.  Suppose $i_j^{\pi_L}<\infty$.
If $b_j<h\le n$ then $v_L(c_h)=nv_K(c_h)\ge n(a_j+1)$,
and hence
\[v_L(c_h)\ge na_j+n>na_j-b_j+h=i_j^{\pi_L}+h.\]
If $1\le h<b_j$ then
\[v_L(c_h)\ge na_j>na_j-b_j+h=i_j^{\pi_L}+h.\]
Finally, we observe that
$v_L(c_{b_j})=na_j=i_j^{\pi_L}+b_j$. \qed \medskip

     For a partition $\lambdabold=
\{\lambda_1,\dots,\lambda_k\}$ whose parts satisfy
$\lambda_i\le n$ for $1\le i\le k$ define
$c_{\lambdabold}=c_{\lambda_1}c_{\lambda_2}\dots
c_{\lambda_k}$.

\begin{prop} \label{Sigma}
Let  $w\ge1$ and let
$\lambdabold=\{\lambda_1,\dots,\lambda_k\}$ be a
partition of $w$ whose parts satisfy $\lambda_i\le n$.
Choose $q$ to minimize $\vbar_p(\lambda_q)$ and
set $t=\vbar_p(\lambda_q)$.  Then 
$v_L(c_{\lambdabold})\ge i_t^{\pi_L}+w$.  If
$v_L(c_{\lambdabold})=i_t^{\pi_L}+w$ and
$i_t^{\pi_L}<\infty$ then $\lambda_q=b_t$ and
$\lambda_i=b_{\nu}=n$ for all $i\not=q$.
\end{prop}

\proof If $i_t^{\pi_L}=\infty$ then
$v_L(c_{\lambda_q})=\infty$, and hence
$v_L(c_{\lambdabold})=\infty$.  Suppose
$i_t^{\pi_L}<\infty$.  By Lemma~\ref{bound} we have
$v_L(c_{\lambda_q})\ge i_t^{\pi_L}+\lambda_q$, and
$v_L(c_{\lambda_i})\ge\lambda_i$ for $i\not=q$.
Hence $v_L(c_{\lambdabold})\ge i_t^{\pi_L}+w$, with
equality if and only if
$v_L(c_{\lambda_q})=i_t^{\pi_L}+\lambda_q$ and
$v_L(c_{\lambda_i})=\lambda_i$ for $i\not=q$.  It
follows from Lemma~\ref{bound} that these conditions
hold if and only if $\lambda_q=b_t$ and
$\lambda_i=b_{\nu}$ for all $i\not=q$. \qed

\begin{prop} \label{ijw}
Let $w\ge1$, let $\mubold$ be a partition of $w$ with
$h\le n$ parts, and set $j=\vbar_p(h)$.  Let
$\lambdabold=\{\lambda_1,\dots,\lambda_k\}$ be a
partition of $w$ whose parts satisfy $\lambda_i\le n$,
choose $q$ to minimize $\vbar_p(\lambda_q)$, and set
$t=\vbar_p(\lambda_q)$.  Then
\\[\medskipamount]
(a) $v_L(d_{\lambdabold\mubold}c_{\lambdabold})
\ge i_j+w$.
\\[\medskipamount]
(b) Suppose
$v_L(d_{\lambdabold\mubold}c_{\lambdabold})=i_j+w$.
Then $i_t^{\pi_L}$ is finite, $\lambda_q=b_t$, and
$\lambda_i=n$ for all $i\not=q$.
\end{prop}

\proof (a) Suppose $t\ge j$.  Then by
Corollary~\ref{val} we have
$\vbar_p(d_{\lambdabold\mubold})\ge t-j$.  Hence by
Proposition~\ref{Sigma} we get
\[v_L(d_{\lambdabold\mubold}c_{\lambdabold})\ge
(t-j)e_L+i_t^{\pi_L}+w\ge i_j+w.\]
Suppose $t<j$.  Using Proposition~\ref{Sigma} we get
\[v_L(d_{\lambdabold\mubold}c_{\lambdabold})\ge
v_L(c_{\lambdabold})\ge i_t^{\pi_L}+w\ge i_t+w
\ge i_j+w.\]
(b) If $v_L(d_{\lambdabold\mubold}c_{\lambdabold})=i_j+w$
then all the inequalities above are equalities.  In
either case it follows that $i_t^{\pi_L}$ is finite and
$v_L(c_{\lambdabold})=i_t^{\pi_L}+w$.  Therefore by
Proposition~\ref{Sigma} we get $\lambda_q=b_t$ and
$\lambda_i=n$ for all $i\not=q$. \qed \medskip

     We now apply some of the results of
section~\ref{digraphs} to our field extension $L/K$.
For $\alpha\in L$ let $M_{\mubold}(\alpha)
=m_{\mubold}(\sigma_1(\alpha),\dots,\sigma_n(\alpha))$.

\begin{prop} \label{spec}
Let $r\ge1$ and let $\alpha\in\M_L^r$.  Choose a power
series
\[\phi(X)=a_rX^r+a_{r+1}X^{r+1}+\dots\]
with coefficients in $\OO_K$ such that
$\alpha=\phi(\pi_L)$.  Then
\[E_h(\alpha)=\sum_{\mubold}
a_{\mu_1}a_{\mu_2}\dots a_{\mu_h}M_{\mubold}(\pi_L),\]
where the sum ranges over all partitions
$\mubold=\{\mu_1,\dots,\mu_h\}$ with $h$ parts such that
$\mu_i\ge r$ for $1\le i\le h$.
\end{prop}

\proof This follows from (\ref{eh}) by setting
$X_i=\sigma_i(\pi_L)$ and letting $a_j\in\OO_K$. \qed

\begin{prop} \label{setting}
Let $n\ge1$, let $w\ge1$, and let $\mubold$ be a
partition of $w$ with at most $n$ parts.  Then
\[M_{\mu}(\pi_L)=\sum_{\lambdabold}
d_{\lambdabold\mubold}c_{\lambdabold},\]
where the sum is over all partitions
$\lambdabold=\{\lambda_1,\dots,\lambda_k\}$ of $w$ such
that $\lambda_i\le n$ for $1\le i\le k$.
\end{prop}

\proof This follows from Theorem~\ref{krtheorem} by
setting $X_i=E_i(\pi_L)=c_i$. \qed \medskip

     Let $1\le h\le n$ and recall that we defined
$g_h:\Z\ra\Z$ by setting $g_h(r)=s$, where $s$ is the
largest integer such that $E_h(\M_L^r)\subset\M_K^s$.

\begin{theorem} \label{contain}
Let $L/K$ be a totally ramified extension of degree
$n=up^{\nu}$, with $p\nmid u$.  Let $r\in\Z$, let
$1\le h\le n$, and set $j=\vbar_p(h)$.  Then
\begin{align*}
E_h(\M_L^r)&\subset\M_K^{\lceil(i_j+hr)/n\rceil} \\
g_h(r)&\ge\left\lceil\frac{i_j+hr}{n}\right\rceil.
\end{align*}
\end{theorem}

\proof Let $\pi_K$ be a uniformizer for $K$.  Then for
$t\in\Z$ we have
\begin{align} \label{Eh}
E_h(\M_L^{nt+r})&=E_h(\pi_K^t\cdot\M_L^r)
=\pi_K^{ht}\cdot E_h(\M_L^r) \\
\left\lceil\frac{i_j+h(nt+r)}{n}\right\rceil
&=ht+\left\lceil\frac{i_j+hr}{n}\right\rceil.
\label{ceil}
\end{align}
Therefore it suffices to prove the theorem in the cases
with $1\le r\le n$.  By Proposition~\ref{spec} each
element of $E_h(\M_L^r)$ is an $\OO_K$-linear
combination of terms of the form $M_{\mubold}(\pi_L)$,
where $\mubold$ is a partition with $h$ parts, all
$\ge r$.  Fix one such partition $\mubold$ and set
$w=\Sigma(\mubold)$; then $w\ge hr$.  Using
Proposition~\ref{setting} we can express
$M_{\mubold}(\pi_L)$ as a sum of terms
$d_{\lambdabold\mubold}c_{\lambdabold}$, where
$\lambdabold=\{\lambda_1,\lambda_2,\dots,\lambda_k\}$ is
a partition of $w$ into parts which are $\le n$.  By
Proposition~\ref{ijw}(a) we get
$v_L(d_{\lambdabold\mubold}c_{\lambdabold})\ge
i_j+w\ge i_j+hr$.  Since
$d_{\lambdabold\mubold}c_{\lambdabold}\in K$ it follows
that $v_K(d_{\lambdabold\mubold}c_{\lambdabold})
\ge\lceil(i_j+hr)/n\rceil$.  Therefore we have
$v_K(M_{\mubold}(\pi_L))\ge \lceil(i_j+hr)/n\rceil$, and
hence $E_h(\M_L^r)\subset\M_K^{\lceil(i_j+hr)/n\rceil}$.
\qed
 
\section{Equality} \label{equal}

In this section we show that in some special cases we
have $\OO_K\cdot E_h(\M_L^r)
=\M_K^{\lceil(i_j+hr)/n\rceil}$, where $j=\vbar_p(h)$.
This is equivalent to showing that
$g_h(r)=\lceil(i_j+hr)/n\rceil$ holds in these cases.
In particular, we prove that if the residue field
$\Kbar$ of $K$ is large enough then
$g_{p^j}(r)=\lceil(i_j+rp^j)/n\rceil$ for
$0\le j\le\nu$.  To prove that
$g_h(r)=\lceil(i_j+hr)/n\rceil$ holds for all $r\in\Z$,
by Theorem~\ref{contain} it suffices to show the
following: Let $r$ satisfy
\begin{equation} \label{break}
\left\lceil\frac{i_j+hr}{n}\right\rceil
<\left\lceil\frac{i_j+h(r+1)}{n}\right\rceil.
\end{equation}
Then there is $\alpha\in\M_L^r$ such that
$v_K(E_h(\alpha))=\lceil(i_j+hr)/n\rceil$.  By
(\ref{Eh}) and (\ref{ceil}) it's enough to prove this
for $r$ such that $1\le r\le n$.

     Once again we let $\pi_L$ be a uniformizer for $L$
whose minimum polynomial over $K$ is
\[f(X)=X^n-c_1X^{n-1}+\cdots+(-1)^{n-1}c_{n-1}X+(-1)^nc_n.\]

\begin{theorem} \label{distinct}
Let $L/K$ be a totally ramified extension of degree
$n=up^{\nu}$, with $p\nmid u$.  Let $j$ be an integer
such that $0\le j\le\nu$ and $\vbar_p(i_j)\ge j$.  Then
\begin{align*}
\OO_K\cdot
E_{p^j}(\M_L^r)&=\M_K^{\lceil(i_j+rp^j)/n\rceil} \\
g_{p^j}(r)&=\left\lceil\frac{i_j+rp^j}{n}\right\rceil.
\end{align*}
\end{theorem}

\proof Set $m=\vbar_p(i_j)$.  Then
$i_j=(m-j)e_L+i_m^{\pi_L}$.  In particular, if
$\Char(K)=p$ then $m=j$ and $i_j=i_m=i_m^{\pi_L}$.
We can write $i_m^{\pi_L}=an-b$ with $1\le b\le n$ and
$\vbar_p(b)=m$.  Since $j\le m$ there is $b'\in\Z$ such
that $b=b'p^j$.  Let $r_1\in\Z$ and set
$r=b'+r_1up^{\nu-j}$.  Then
\begin{equation} \label{ijrpj}
i_j+rp^j=(m-j)e_L+an+r_1n.
\end{equation}
Therefore we have
\begin{align*}
\left\lceil\frac{i_j+rp^j}{n}\right\rceil
&=(m-j)e_K+a+r_1 \\
\left\lceil\frac{i_j+(r+1)p^j}{n}\right\rceil
&=(m-j)e_K+a+r_1+1,
\end{align*}
with $e_K=v_K(p)=e_L/n$.  It follows that the only
values of $r$ in the range $1\le r\le n$ satisfying
(\ref{break}) are of the form $r=b'+r_1up^{\nu-j}$ with
$0\le r_1<p^j$.  Therefore it suffices to prove that
$v_K(E_{p^j}(\pi_L^r))=(m-j)e_K+a+r_1$ holds for these
values of $r$.

     Let $\mubold$ be the partition of $rp^j$ consisting
of $p^j$ copies of $r$.  Then $E_{p^j}(\pi_L^r)
=M_{\mubold}(\pi_L)$, so it follows from
Proposition~\ref{setting} that
\begin{equation} \label{Epjsum0}
E_{p^j}(\pi_L^r)=\sum_{\lambdabold}
d_{\lambdabold\mubold}c_{\lambdabold},
\end{equation}
where the sum is over all partitions
$\lambdabold=\{\lambda_1,\dots,\lambda_k\}$ of $rp^j$
such that $\lambda_i\le n$ for $1\le i\le k$.
It follows from Proposition~\ref{ijw}(a) that
$v_L(d_{\lambdabold\mubold}c_{\lambdabold})\ge i_j+rp^j$.
Suppose $v_L(d_{\lambdabold\mubold}c_{\lambdabold})
=i_j+rp^j$.  Then by Proposition~\ref{ijw}(b) we see
that $\lambdabold$ has at most one element which is not
equal to $n$.  Since $\Sigma(\lambdabold)=rp^j=b+r_1n$,
and the elements of $\lambdabold$ are $\le n$, it
follows that $\lambdabold=\kappabold$, where
$\kappabold$ is the partition of $rp^j$ which consists
of 1 copy of $b$ and $r_1$ copies of $n$.  Since
$E_{p^j}(\pi_L^r)\in K$ and $d_{\kappabold\mubold}
c_{\kappabold}\in K$ it follows from (\ref{Epjsum0}) and
(\ref{ijrpj}) that
\begin{equation} \label{Epj0}
E_{p^j}(\pi_L^r)\equiv d_{\kappabold\mubold}
c_{\kappabold}\pmod{\M_K^{(m-j)e_K+a+r_1+1}}.
\end{equation}

     Let $\kappabold'$ be the partition of $r$
consisting of 1 copy of $b'$ and $r_1$ copies of
$up^{\nu-j}$, and let $\mubold'$ be the partition of $r$
consisting of 1 copy of $r$.  Then
$\kappabold=p^j\cdot\kappabold'$ and
$\mubold=p^j*\mubold'$.  Since $v_p(b')=m-j$ it follows
from Proposition~\ref{cong} that
$d_{\kappabold\mubold}\equiv d_{\kappabold'\mubold'}
\pmod{p^{m-j+1}}$.  Suppose $m<\nu$.  Then $b<n$, so
$b'\not=up^{\nu-j}$.  Hence by Proposition~\ref{cycle1}
we get $d_{\kappabold'\mubold'}=(-1)^{r+r_1+1}r$.  Since
$r=b'+r_1up^{\nu-j}$ and $v_p(b')=m-j$ this implies
$v_p(d_{\kappabold'\mubold'})=v_p(r)=m-j$.  Suppose
$m=\nu$.  Then $b=n$ and $b'=p^{-j}b=up^{\nu-j}$, so
$\kappabold'$ consists of $r_1+1$ copies of
$up^{\nu-j}$.  Since $\gcd(up^{\nu-j},r)=up^{\nu-j}$
and $\gcd(r_1+1,1)=1$, by Proposition~\ref{cycle3} we
get $d_{\kappabold'\mubold'}=(-1)^{r+r_1+1}up^{\nu-j}$.
Hence $v_p(d_{\kappabold'\mubold'})=\nu-j=m-j$ holds in
this case as well.  Since $d_{\kappabold\mubold}\equiv
d_{\kappabold'\mubold'}\pmod{p^{m-j+1}}$ it follows that
$v_p(d_{\kappabold'\mubold'})=m-j$.  Therefore
$v_K(d_{\kappabold\mubold}c_{\kappabold})=(m-j)e_K+a+r_1$.
Using (\ref{Epj0}) we conclude that
$v_K(E_{p^j}(\pi_L^r))=(m-j)e_K+a+r_1$.~\qed

\begin{theorem} \label{multiple}
Let $L/K$ be a totally ramified extension of degree
$n=up^{\nu}$, with $p\nmid u$.  Let $j$ be an integer
such that $0\le j\le\nu$ and $\vbar_p(i_j)<j$.  Set
$m=\vbar_p(i_j)$ and assume that $|\Kbar|>p^m$.  Then
\begin{align*}
\OO_K\cdot E_{p^j}(\M_L^r)
&=\M_K^{\lceil(i_j+rp^j)/n\rceil} \\
g_{p^j}(r)&=\left\lceil\frac{i_j+rp^j}{n}\right\rceil.
\end{align*}
\end{theorem}

\proof Since $m<j$ we have $i_m=i_j=i_j^{\pi_L}$.
Therefore $i_j=an-b$ for some $a,b$ such that $1\le b<n$
and $\vbar_p(b)=m$.  Hence $b=b'p^j+b''p^m$ for some
$b',b''$ such that $0<b''<p^{j-m}$ and $p\nmid b''$.
Let $r_1\in\Z$ and set $r=b'+r_1up^{\nu-j}$.  Then
\begin{equation} \label{ijrp}
i_j+rp^j=an+r_1n-b''p^m,
\end{equation}
so we have
\begin{align*}
\left\lceil\frac{i_j+rp^j}{n}\right\rceil
&=a+r_1+\left\lceil\frac{-b''p^m}{n}\right\rceil
=a+r_1 \\[.2cm]
\left\lceil\frac{i_j+(r+1)p^j}{n}\right\rceil
&=a+r_1+\left\lceil\frac{p^j-b''p^m}{n}\right\rceil
=a+r_1+1.
\end{align*}
It follows that the only values of $r$ in the range
$1\le r\le n$ satisfying (\ref{break}) are of the form
$r=b'+r_1up^{\nu-j}$ with $0\le r_1<p^j$.  It suffices
to prove that for every such $r$ there is
$\beta\in\OO_K$ such that
$v_K(E_{p^j}(\pi_L^r+\beta\pi_L^{r+b''}))=a+r_1$.

     Let $\eta(X)=E_{p^j}(\pi_L^r+X\pi_L^{r+b''})$.
We need to show that there is $\beta\in\OO_K$ such that
$v_K(\eta(\beta))=a+r_1$.  It follows from
Proposition~\ref{spec} that $\eta(X)$ is a polynomial in
$X$ of degree at most $p^j$, with coefficients in
$\OO_K$.  For $0\le\ell\le p^j$ let
$\mubold^{\ell}$ be the partition of $rp^j+\ell b''$
consisting of $p^j-\ell$ copies of $r$ and $\ell$ copies
of $r+b''$.  By Proposition~\ref{spec} the coefficient
of $X^\ell$ in $\eta(X)$ is equal to
$M_{\mubold^\ell}(\pi_L)$.  By Proposition~\ref{setting}
we have
\begin{equation} \label{Mmul}
M_{\mubold^\ell}(\pi_L)=\sum_{\lambdabold}
d_{\lambdabold\mubold^{\ell}}c_{\lambdabold},
\end{equation}
where the sum is over all partitions
$\lambdabold=\{\lambda_1,\dots,\lambda_k\}$ of
$rp^j+\ell b''$ such that $\lambda_i\le n$ for
$1\le i\le k$.  Using Proposition~\ref{ijw}(a) and
(\ref{ijrp}) we get
\begin{align} \nonumber
v_L(d_{\lambdabold\mubold^{\ell}}c_{\lambdabold})
&\ge i_j+rp^j+\ell b'' \\
&=(a+r_1)n+(\ell-p^m)b'' \label{first} \\
&>(a+r_1-1)n. \nonumber
\end{align}
Since $d_{\lambdabold\mubold^{\ell}}c_{\lambdabold}
\in K$ it follows that $d_{\lambdabold\mubold^{\ell}}
c_{\lambdabold}\in\M_K^{a+r_1}$.  Therefore by
Proposition~\ref{setting} we have
$M_{\mubold^\ell}(\pi_L)\in\M_K^{a+r_1}$.

     Suppose $v_K(d_{\lambdabold\mubold^{\ell}}
c_{\lambdabold})=a+r_1$.  Then
$v_L(d_{\lambdabold\mubold^{\ell}}c_{\lambdabold})
=(a+r_1)n$, so by (\ref{first}) we get $\ell\le p^m$.
Hence for $p^m<\ell\le p^j$ we have
$M_{\mubold^\ell}(\pi_L)\in\M_K^{a+r_1+1}$.  Let
$w=b+r_1n=rp^j+b''p^m$ and let $\mubold=\mubold^{p^m}$
be the partition of $w$ consisting of $p^m$ copies of
$r+b''$ and $p^j-p^m$ copies of $r$.  Then the
coefficient of $X^{p^m}$ in $\eta(X)$ is
$M_{\mubold}(\pi_L)$.  Let $\kappabold$ be the partition
of $w$ consisting of 1 copy of $b$ and $r_1$ copies of
$n$.  Suppose $\lambdabold$ is a partition of $w$ with
parts $\le n$ such that
$v_K(d_{\lambdabold\mubold}c_{\lambdabold})=a+r_1$.
Since $(a+r_1)n=i_j+w$ it follows from
Proposition~\ref{ijw}(b) that $\lambdabold$ has
at most one element which is not equal to $n$.  Since
$\Sigma(\lambdabold)=b+r_1n$, and the elements of
$\lambdabold$ are $\le n$, it follows that
$\lambdabold=\kappabold$.  Hence by (\ref{Mmul}) we have
\begin{equation} \label{Mcong}
M_{\mubold}(\pi_L)\equiv d_{\kappabold\mubold}
c_{\kappabold}\pmod{\M_K^{a+r_1+1}}.
\end{equation}

     Set $w'=b'p^{j-m}+b''+r_1up^{\nu-m}=rp^{j-m}+b''$.
Let $\kappabold'$ be the partition of $w'$ consisting of
1 copy of $b'p^{j-m}+b''$ and $r_1$ copies of
$up^{\nu-m}$, and let $\mubold'$ be the partition of
$w'$ consisting of 1 copy of $r+b''$ and $p^{j-m}-1$
copies of $r$.  Then $\kappabold=p^m\cdot\kappabold'$
and $\mubold=p^m*\mubold'$, so by Proposition~\ref{cong}
we have $d_{\kappabold\mubold}\equiv
d_{\kappabold'\mubold'}\pmod{p}$.
Let $\Gamma$ be a cycle digraph which has an
admissible $(\kappabold',\mubold')$-tiling.
Suppose $\Gamma$ has more than one component.  Since
$\Gamma$ has a $\kappabold$-tiling, $\Gamma$ has at
least one component $\Gamma_0$ such that $up^{\nu-m}$
divides $|V(\Gamma_0)|$.  Thus
$|V(\Gamma_0)|=k\cdot up^{\nu-m}$ for some $k$ such that
$1\le k\le r_1$.  Let $\kappabold_0'$ be the submultiset
of $\kappabold'$ consisting of $k$ copies of
$up^{\nu-m}$.  Then $\kappabold_0'$ is the unique
submultiset of $\kappabold'$ such that $\Gamma_0$ has a
$\kappabold_0'$-tiling.  Furthermore there is a unique
submultiset $\mubold_0'$ of $\mubold'$ such that
$\Gamma_0$ has a $\mubold_0'$-tiling.

     Suppose $r$ does not divide $kup^{\nu-m}$.  Then
there is $\ell\ge0$ such that $\mubold_0'$ consists of 1
copy of $r+b''$ together with $\ell$ copies of $r$.  By
Proposition~\ref{cycle0} we have
$\eta_{\kappabold_0'\mubold_0'}(\Gamma_0)=up^{\nu-m}$.
Let $\Gamma_1$ be the complement of $\Gamma_0$ in
$\Gamma$, let $\kappabold_1'=\kappabold'\smallsetminus
\kappabold_0'$, and let $\mubold_1'=\mubold'
\smallsetminus\mubold_0'$.  Since $\Gamma_1$ has no
cycle of length $|V(\Gamma_0)|=b''+(\ell+1)r$
we have $\eta_{\kappabold'\mubold'}(\Gamma)
=\eta_{\kappabold_0'\mubold_0'}(\Gamma_0)
\eta_{\kappabold_1'\mubold_1'}(\Gamma_1)$.
Hence $\eta_{\kappabold'\mubold'}(\Gamma)$ is
divisible by $p$ in this case.

     On the other hand, suppose $r$ divides
$kup^{\nu-m}$.  Then there is $\ell\ge1$ such that
$\mubold_0'$ consists of $\ell$ copies of $r$.  It
follows that
$k\cdot up^{\nu-m}=\ell\cdot r$.  Let $(S,T)$ be an
admissible $(\kappabold',\mubold')$-tiling of $\Gamma$
and let $(S_0,T_0)$ be the restriction of $(S,T)$ to
$\Gamma_0$.  Then $(S_0,T_0)$ is a
$(\kappabold_0',\mubold_0')$-tiling of $\Gamma_0$.
By Proposition~\ref{cycle2}(b) the
automorphism group of $(\Gamma_0,S_0,T_0)$ has order
$\gcd(k,\ell)$.  Since $\Aut(\Gamma_0,S_0,T_0)$ is
isomorphic to a subgroup of $\Aut(\Gamma,S,T)$, it
follows that $\gcd(k,\ell)$ divides
$|\Aut(\Gamma,S,T)|$.  Therefore the assumption that
$(S,T)$ is admissible implies
that $\gcd(k,\ell)=1$.  Hence $k\mid r$ and $\ell\mid
up^{\nu-m}$, so there is $q\in\Z$ with $r=kq$ and
$up^{\nu-m}=\ell q$.  By Proposition~\ref{cycle2}(a) the
number of isomorphism classes of
$(\kappabold_0',\mubold_0')$-tilings of $\Gamma_0$ is
\[\eta_{\kappabold_0'\mubold_0'}(\Gamma_0)=
\gcd(up^{\nu-m},r)=\gcd(\ell q,kq)=q.\]
If $p\mid q$ then as above we deduce that
$\eta_{\kappabold'\mubold'}(\Gamma)$ is divisible by
$p$.  On the other hand, if $p\nmid q$ then $q\mid u$;
in particular, $q\le u$.  Since $k\le r_1$ we get
$r_1up^{\nu-j}+b'=r=kq\le r_1u$, a contradiction.  By
combining the two cases we find that if $\Gamma$ has
more than one component then
$\eta_{\kappabold'\mubold'}(\Gamma)$ is divisible by
$p$.

     Finally, suppose that $\Gamma$ consists of a single
cycle of length $w'$.  Then by Proposition~\ref{cycle4}
we have $\eta_{\kappabold',\mubold'}(\Gamma)=w'$.
Hence by (\ref{dlm}) we get
\[d_{\kappabold\mubold}\equiv d_{\kappabold'\mubold'}
\equiv\pm\eta_{\kappabold'\mubold'}(\Gamma)
\equiv\pm w'\pmod{p}.\]
Since $w'\equiv b''\pmod{p}$ it follows that
$p\nmid d_{\kappabold\mubold}$.  Hence by (\ref{Mcong})
we get
\[v_K(M_{\mubold}(\pi_L))=v_K(c_{\kappabold})=a+r_1.\]
Let $\pi_K$ be a uniformizer for $K$ and set
$\phi(X)=\pi_K^{-a-r_1}\eta(X)$.  Then
$\phi(X)\in\OO_K[X]$.  Let $\phibar(X)$ be the image of
$\phi(X)$ in $\Kbar[X]$.  We have shown that
$\phibar(X)$ has degree $p^m$.  Since $|\Kbar|>p^m$
there is $\betabar\in\Kbar$ such that
$\phibar(\betabar)\not=0$.  Let $\beta\in\OO_K$ be a
lifting of $\betabar$.  Then
$\phi(\beta)\in\OO_K^{\times}$.  It follows that
\[v_K(E_{p^j}(\pi_L^r+\beta\pi_L^{r+b''}))
=v_K(\eta(\beta))=a+r_1.\]
We conclude that $\OO_K\cdot E_{p^j}(\M_L^r)
=\M_K^{\lceil(i_j+rp^j)/n\rceil}$. \qed \medskip

\begin{remark}
Theorems~\ref{distinct} and \ref{multiple} together
imply that if $\Kbar$ is sufficiently large then
$g_{p^j}(r)=\lceil(i_j+rp^j)/n\rceil$ for
$0\le j\le\nu$.  This holds for instance if
$|\Kbar|\ge p^{\nu}$.
\end{remark} \medskip

\begin{remark} \label{dj}
Let $L/K$ be a totally ramified separable extension of
degree $n=up^{\nu}$.  The different $\M_L^{d_0}$ of
$L/K$ is defined by letting $d_0$ be the largest integer
such that $E_1(\M_L^{-d_0})\subset\OO_K$.  For
$1\le j\le\nu$ one can define higher order analogs
$\M_L^{d_j}$ of the different by letting $d_j$ be the
largest integer such that
$E_{p^j}(\M_L^{-d_j})\subset\OO_K$.  An argument similar
to the proof of Proposition~\ref{different} shows that
\[\OO_K\cdot E_{p^j}(\M_L^r)
=\M_K^{\lfloor p^j(d_j+r)/n\rfloor}.\]
This generalizes Proposition~\ref{different}, which
is equivalent to the case $j=0$ of this formula.  By
Proposition~3.18 of \cite{heier}, the valuation of the
different of $L/K$ is $d_0=i_0+n-1$.
Using Theorems~\ref{distinct} and \ref{multiple} we find
that, if $\Kbar$ is sufficiently large, $d_j$ is the
largest integer such that
${\lceil(i_j-d_jp^j)/n\rceil}\ge0$.  Hence
$d_j=\lfloor(i_j+n-1)/p^j\rfloor$ for $0\le j\le\nu$.
\end{remark} \medskip

\begin{example}
Let $K=\F_2((t))$ and let $L$ be an extension of $K$
generated by a root $\pi_L$ of the Eisenstein polynomial
$f(X)=X^8+tX^3+tX^2+t$.  Then the indices of
inseparability of $L/K$ are $i_0=3$, $i_1=i_2=2$, and
$i_3=0$.  Since $\lceil(i_2+2^2\cdot1)/2^3\rceil=1$, the
formula in Theorem~\ref{multiple} would imply
$\OO_K\cdot E_4(\M_L^1)=\M_K^1$.  We claim that
$E_4(\M_L)\subset\M_K^2$.

     Let $\alpha\in\M_L$ and write
$\alpha=a_1\pi_L+a_2\pi_L^2+\dots$, with $a_i\in\F_2$.
It follows from Propositions~\ref{spec} and
\ref{setting} that $E_4(\alpha)$ is a sum of terms of
the form
$a_{\mubold}d_{\lambdabold\mubold}c_{\lambdabold}$,
where $\lambdabold$ is a partition whose parts are
$\le8$ and $\mubold$ is a partition with 4 parts such
that $\Sigma(\lambdabold)=\Sigma(\mubold)$.  We are
interested only in those terms with $K$-valuation 1.  We
have $v_K(c_{\lambdabold})\ge2$ unless $\lambdabold$ is
one of $\{5\}$, $\{6\}$, or $\{8\}$.  If
$\lambdabold=\{8\}$
then $2\mid d_{\lambdabold\mubold}$ for any $\mubold$ by
Corollary~\ref{val}.  If $\lambdabold=\{6\}$ and
$\mubold=\{1,1,1,3\}$ then $d_{\lambdabold\mubold}=6$ by
Proposition~\ref{cycle1}.  If $\lambdabold=\{6\}$ and
$\mubold=\{1,1,2,2\}$ then a computation based on
(\ref{dlm}) shows that $d_{\lambdabold\mubold}=9$.
If $\lambdabold=\{5\}$ and $\mubold=\{1,1,1,2\}$ then
$d_{\lambdabold\mubold}=-5$ by Proposition~\ref{cycle1}.
Combining these facts we get
\[E_4(\alpha)\equiv a_1^3a_2t+a_1^2a_2^2t
\pmod{\M_K^2}.\]
Since $a_1,a_2\in\F_2$ we have $a_1^3a_2+a_1^2a_2^2=0$.
Therefore $E_4(\alpha)\in\M_K^2$.  Since this holds for
every $\alpha\in\M_L$ we get $E_4(\M_L)\subset\M_K^2$.
This shows that Theorem~\ref{multiple} does not hold
without the assumption about the size of $\Kbar$.
\end{example} \medskip

     The following result shows that
$g_h(r)=\lceil(i_j+hr)/n\rceil$ does not hold in
general, even if we assume that the residue field of
$K$ is large.  It also suggests that there may not be a
simple criterion for determining when
$g_h(r)=\lceil(i_j+hr)/n\rceil$ does hold.

\begin{prop}
Let $L/K$ be a totally ramified extension of degree $n$,
with $p\nmid n$.  Let $r\in\Z$ and $1\le h\le n$ be such
that $n\mid hr$.  Set $s=hr/n$, $u=\gcd(r,n)$, and
$v=\gcd(h,s)$.  Then $g_h(r)=\lceil(i_0+hr)/n\rceil=s$
if and only if $p$ does not divide the binomial
coefficient $\dst\binom{u}{v}$.  In particular, if $u<v$
then $g_h(r)>s$.
\end{prop}

\proof Since $L/K$ is tamely ramified we have $\nu=0$,
$i_0=0$, and
\[\left\lceil\frac{i_0+hr}{n}\right\rceil
=\left\lceil\frac{hr}{n}\right\rceil=s.\]
It follows from Theorem~\ref{contain} that
$g_h(r)\ge s$.  If $r'=nt+r$ then $s'=hr'/n=ht+s$,
$u'=\gcd(r',n)=u$, and $v'=\gcd(h,s')=v$.  Hence by
(\ref{Eh}) it suffices to prove the proposition in the
cases with $1\le r\le n$.

     Suppose $p$ does not divide $\dst\binom{u}{v}$.
To prove $g_h(r)=s$ it suffices to show that
$v_K(E_h(\pi_L^r))=s$.  Let $\mubold$ be the partition
of $hr$ consisting of $h$ copies of $r$.  Then
$E_h(\pi_L^r)=M_{\mubold}(\pi_L)$, so it follows from
Proposition~\ref{setting} that
\begin{equation} \label{Ehsum}
E_h(\pi_L^r)=\sum_{\lambdabold}
d_{\lambdabold\mubold}c_{\lambdabold},
\end{equation}
where the sum is over all partitions
$\lambdabold=\{\lambda_1,\dots,\lambda_k\}$ of $hr$
such that $\lambda_i\le n$ for $1\le i\le k$.
Let $\kappabold$ be the partition of $hr=sn$ consisting
of $s$ copies of $n$ and let $\lambdabold$ be a partition
of $hr$ whose parts are $\le n$.  Then by
Proposition~\ref{ijw}(a) we have
$v_L(d_{\kappabold\mubold}c_{\lambdabold})\ge hr=sn$.
Furthermore, if
$v_L(d_{\kappabold\mubold}c_{\lambdabold})=hr$ then by
Proposition~\ref{ijw}(b) we have
$\lambdabold=\kappabold$.  Hence by (\ref{Ehsum}) we get
\[E_h(\pi_L^r)\equiv d_{\kappabold\mubold}c_{\kappabold}
\pmod{\M_K^{s+1}}.\]
By Proposition~\ref{cycle3} we have
$d_{\kappabold\mubold}=\pm\dst\binom{u}{v}$.  Since
$\dst p\nmid\binom{u}{v}$ and $v_K(c_{\kappabold})=s$ it
follows that $v_K(E_h(\pi_L^r))=s$.  Therefore
$g_h(r)=s$.

     Suppose $p$ divides $\dst\binom{u}{v}$.  By
Proposition~\ref{spec}, each element of $E_h(\M_L^r)$
is an $\OO_K$-linear combination of terms of the form
$M_{\nubold}(\pi_L)$ where $\nubold$ is a partition
with $h$ parts, all $\ge r$.  Fix one such partition
$\nubold$ and set $w=\Sigma(\nubold)$; then $w\ge
hr=sn$.  By Proposition~\ref{setting} we can express
$M_{\nubold}(\pi_L)$ as a sum of terms of the form
$d_{\lambdabold\nubold}c_{\lambdabold}$, where
$\lambdabold=\{\lambda_1,\lambda_2,\dots,\lambda_k\}$ is
a partition of $w$ into parts which are $\le n$.  By
Proposition~\ref{ijw}(a) we have
$v_L(d_{\lambdabold\nubold}c_{\lambdabold})\ge w\ge sn$.
Suppose $v_L(d_{\lambdabold\nubold}c_{\lambdabold})=sn$.
Then $w=sn$, and by Proposition~\ref{ijw}(b) we see that
$\lambdabold$ consists of $k$ copies of $n$.
It follows that $kn=w=sn$, and hence that $k=s$.
Therefore $\lambdabold=\kappabold$.  Since
$\Sigma(\nubold)=w=kn=hr$ we get
$\nubold=\mubold$.  Since
$d_{\kappabold\mubold}=\pm\dst\binom{u}{v}$ and $p$
divides $\dst\binom{u}{v}$ we have
$v_L(d_{\kappabold\mubold}c_{\kappabold})>
v_L(c_{\kappabold})=sn$, a
contradiction.  Hence
$v_L(d_{\lambdabold\nubold}c_{\lambdabold})>sn$ holds in
all cases.  Since $d_{\lambdabold\nubold}c_{\lambdabold}
\in K$ we get $v_K(d_{\lambdabold\nubold}c_{\lambdabold})
\ge s+1$.  It follows that
$E_h(\M_L^r)\subset\M_K^{s+1}$, and hence that
$g_h(r)\ge s+1$. \qed


\begin{thebibliography}{9}

\bibitem{FV} I. B. Fesenko and S. V. Vostokov,
{\em Local fields and their extensions.  A constructive
approach}, Amer.\ Math.\ Soc., Providence, RI, 2002.

\bibitem{fried} M. Fried, Arithmetical properties of
function fields II, The generalized Schur problem, Acta
Arith.\ {\bf25} (1973/74), 225--258.

\bibitem{fm} M. Fried and A. M\'ezard, Configuration
spaces for wildly ramified covers, appearing in
{\em Arithmetic Fundamental Groups and Noncommutative
Algebra}, Proc.\ Sympos.\ Pure Math.\ {\bf70} (2002),
353--376.

\bibitem{heier} V. Heiermann, De nouveaux invariants
num\'eriques pour les extensions totalement ramifi\'ees
de corps locaux, J. Number Theory {\bf 59} (1996),
159--202.

\bibitem{towers} K. Keating, Indices of inseparability
in towers of field extensions, J. Number Theory {\bf150}
(2015), 81--97.

\bibitem{kr} A. Kulikauskas and J. Remmel, Lyndon words
and transition matrices between elementary, homogeneous
and monomial symmetric functions, Electronic J.
Combinatorics {\bf13} (2006), \#R18.

\bibitem{cl} J.-P. Serre, {\em Corps Locaux}, Hermann,
Paris (1962).

\end{thebibliography}
\end{document}